\documentclass{elsart}

\usepackage{amsmath}
\usepackage{amssymb}
\usepackage{graphicx}
\usepackage{color}

\newenvironment{biography}[2]{
\footnotesize \unitlength 1mm 
\bigskip \parskip=0pt \par 
\rule{0pt}{39mm} \vspace{-39mm} \par
\noindent\setbox0\hbox{{#1}}
\ht0=39mm \count10=\ht0 \divide \count10 by \baselineskip
\global \hangindent29mm \global \hangafter- \count10%
\hskip-28.5mm \setbox0 \hbox to 28.5mm {\raise-30.5mm\box0\hss}%
\dp0=0mm \ht0=0mm \box0 \noindent \bf#2 \rm}
{\par\rm\small} 

\begin{document}

\begin{frontmatter}

\journal{Neurocomputing}

\title{Forecasting the CATS benchmark with the Double Vector Quantization method}

\author[GS]{Geoffroy Simon\corauthref{GS-cor}\thanksref{GS-thk}} 
\ead{simon@dice.ucl.ac.be}
\ead[url]{http://www.dice.ucl.ac.be/\~{}gsimon/}
\thanks[GS-thk]{G. Simon is funded by the Belgian F.R.I.A.}
\corauth[GS-cor]{Corresponding author: G. Simon}
\author[GS]{John A. Lee\thanksref{JL-thk}} 
\thanks[JL-thk]{J. A. Lee is a Scientific Research Worker of the Belgian F.N.R.S.}
\author[MC]{Marie Cottrell} 
\author[GS,MC]{Michel Verleysen\thanksref{MV-thk}}
\thanks[MV-thk]{M. Verleysen is a Research Director of the Belgian F.N.R.S.}

\address[GS]{Universit\'e catholique de Louvain, Machine Learning Group - DICE\\ 
	Place du Levant 3, B-1348 Louvain-la-Neuve, Belgium}
\address[MC]{Universit\'e Paris I - Panth\'eon Sorbonne, SAMOS-MATISSE, UMR CNRS 8595\\
	Rue de Tolbiac 90, F-75634 Paris Cedex 13, France}

\begin{abstract}
The Double Vector Quantization method, a long-term forecasting method based 
on the SOM algorithm, has been used to predict the 100 missing values of the CATS 
competition data set. An analysis of the proposed time series is provided 
to estimate the dimension of the auto-regressive part of this nonlinear 
auto-regressive forecasting method. Based on this analysis experimental results 
using the Double Vector Quantization (DVQ) method are presented and discussed. 
As one of the features of the DVQ method is its ability to predict scalars 
as well as vectors of values, the number of iterative predictions needed 
to reach the prediction horizon is further observed. The method stability 
for the long term allows obtaining reliable values for a rather long-term 
forecasting horizon.
\end{abstract}


\end{frontmatter}

\section{Introduction}\label{intro}
Time series prediction is a problem encountered in many fields: from engineering 
(predictive control of industrial processes) to finance (forecasting returns of shares 
or stock markets) this general problem has already been studied by a large community
of researchers. Models and prediction methodologies have been proposed by statisticians, 
mathematicians and engineers, as well as people from econometrics and more recently 
from the neural network community. Whatever the time series, whatever the method used 
to predict the series, the methodology always first consists in constructing a model 
of the time series. This model is then used to predict the future of the series. 
In this paper, this general methodology will be applied, using a nonlinear auto-regressive 
model, namely the Double Vector Quantization (DVQ) method \cite{Simon04}.

This paper concerns the prediction of the CATS benchmark proposed as a time series competition 
at IJCNN 2004. As explained in the descriptive paper \cite{Lendasse04}, the CATS competition 
prediction problem consists in predicting 100 missing values distributed in five gaps 
of twenty data points within a time series of 4900 data points. Each gap is preceded 
by 980 known values. The data set starts with 980 known data and the last gap is 
at the end of the time series. In this paper, a global model able to predict the 100 
missing values is developed. To set the parameters of the model (number of prototypes 
in the quantization steps, choice of the preprocessing and of the block prediction horizon), 
an extended cross-validation procedure is used. 

The DVQ nonlinear auto-regressive \cite{Simon04} method is based on Kohonen's self-organizing 
maps (SOM) \cite{Kohonen95}. Two SOMs are used and linked by a stochastic model. The two SOMs 
are used here as a clustering tool for prediction even if the SOM algorithm is usually 
considered as a classification or feature extraction tool. Nevertheless some previous 
attempts have already tried to use the SOM algorithm for time series prediction. For example 
\cite{Cottrell96} uses SOMs to create clusters in the regressor space, \cite{Walter90} and 
\cite{Vesanto97} associate each cluster to a linear local model and \cite{Dablemont03} associates 
each cluster to a nonlinear one. Another approach is to split the problem into the predictions 
of respectively a normalized curve and the curve mean and standard deviation \cite{Cottrell}. 
Using the recursive SOMs approach \cite{Voegtlin02} (and pioneer work on leaky integrators 
\cite{Chappell93}) one tries to learn sequences of data, as applied in \cite{Kangas94} 
for speech recognition problems. Recursive SOMs can be further combined with local linear 
models, as in \cite{Koskela98}. The DVQ method is rather different from these previous works 
since it aims at predicting long term trends instead of providing short term accurate 
predictions. Furthermore, as long term predictions are the main concern, a theoretical proof 
of the DVQ stability at long term has been proposed recently \cite{Simon04}. Considering 
that a gap of 20 values is a long-term horizon when compared to the usual one-step-ahead 
framework, the DVQ method is used to forecast the 100 missing values of the CATS data set.

Another original aspect of the DVQ method is its ability to predict vectors of values 
in a single step. This generalization from scalar to vector cases will be used and studied 
on the CATS data set: the DVQ method can recursively predict 20 times a one-step-ahead scalar 
forecast, or a vector of 20 values in one single prediction step, or any other intermediate 
situation. Note that, as a consequence, the size of the prediction vectors is a parameter 
of the DVQ model. This parameter will be optimized by cross-validation in the experiments.

In the following of this paper, we first present an analysis of the CATS data set. In this 
analysis the correlation dimension of the CATS time series is estimated. The goal of this 
analysis is to determine the size of the regressors for the DVQ nonlinear auto-regressive model. 
In section \ref{method} a short reminder on the basic concepts about the SOMs is followed 
by the description of the forecasting method for the scalar case. Its extension to vectors 
is briefly sketched, and some general comments on the DVQ method and its use in practice 
are given. Section \ref{Methodo} is devoted to the description of the experimental methodology 
that has been specifically developed and used for the CATS series. Section \ref{results} 
presents the results, before a concluding discussion.

\section{Analysis of the CATS data set}\label{analyseCATS}
As mentioned in the introduction, the model implemented by the DVQ method is a NAR 
model. As usual with (N)AR models, the key point is to determine the order of the AR part, 
i.e. to evaluate the size of the regressor. More formally, having at disposal a time series 
of $x(t)$ values with $1 \leq t \leq n$, predicting the values for $t > n$ can be defined as:
\begin{equation}\label{pred}
[x(t+1), ..., x(t+d)] = f(x(t), ..., x(t-p+1), \theta) + \varepsilon_t,
\end{equation}
where $d$ is the size of the vector of values to be predicted, $f$ is the model of 
the data generating process (considered here as a nonlinear one), $p$ is the number 
of past values to consider, $\theta$ are the model parameters and $\varepsilon_t$ 
is the noise. The past values are gathered in a $p$-dimensional vector called 
{\it regressor}. Both $p$ and $d$ must be chosen. As mentioned above, $d$ will result 
from a compromise between scalar predictions with repetitions and a larger vector 
of values to predict as a whole; in practice the choice of $d$ will be determined 
by extensive simulations. 

Concerning the value $p$ of the regressor size, it is possible to have insights about 
a plausible value (or range of values) by an in-depth examination of the series. 
In this paper, the search for the regressor size $p$ is based on Grassberger-Proccacia's 
\cite{Grass} correlation dimension; the procedure is for example summarized in \cite{GrassNPL}. 

Grassberger-Proccacia's procedure allows estimating the correlation dimension $D_c$ \cite{Grass} 
of a time series. Then, according to Takens theorem \cite{Takens}, a regressor of size 
$p = 2 * D_c+1$ will describe the data in an embedding space containing enough information 
to allow a modeling of the manifold describing by the time series values.

In short, the Grassberger-Procaccia procedure computes the correlation dimension $D_c$ 
of vector data $x_t$ according to:
\begin{equation}\label{dim_corr}
\displaystyle D_c = \lim_{r \to 0} \frac{ln(C(r))}{ln(r)},
\end{equation}
where $C(r)$ is the {\it correlation integral} \cite{Grass} defined as:
\begin{equation}\label{int_corr}
\displaystyle C(r) = \lim_{n \to \infty} \frac{2}{n(n - 1)} \sum_{1 \leq t < t' \leq n}^{} I\left(\| x_t - x_{t'}\| \leq r\right).
\end{equation}
Function $I(.)$ takes a value equal to 1 if its expression into parenthesis is true 
and 0 otherwise.

Intuitively, the idea in relation (\ref{int_corr}) is to count the number of data $x_{t'}$ 
in a hypersphere centered in $x_t$ with radius $r$. This operation is repeated for each data 
$x_t$. Then the limit for $n$ tends to $\infty$ is taken, i.e. the definition is given 
for an infinite number of data $x_t$ (in the series). Finally, the ratio between the log 
of this number of data and the log of the corresponding radius is observed in relation 
(\ref{dim_corr}), as the radius $r$ tends to zero. In other words while estimating 
the correlation integral (\ref{int_corr}), one tries to count the number of data $x_{t'}$ 
that are at most at distance $r$ from $x_t$, given an infinite number of data. The correlation 
integral thus expresses the asymptotic proportion of data pairs whose distance is less than $r$, 
with repect to the total number $\frac{n(n-1)}{2}$ of pairs. The correlation dimension 
(\ref{dim_corr}) is therefore obtained as the limit of the correlation integral values obtained 
when the radius $r$ is decreased to zero. As mentioned in \cite{GrassNPL} the correlation dimension 
is given by the slope in the linear part of the curves. This follows the initial assumption 
by Grassberger and Procaccia \cite{Grass} that the correlation integral $C(r)$ behaves 
as a power law of $r$ for small values of $r$. In other words, $C(r) \approx r^{D_c}$. 
The correlation dimension is therefore obtained as the slope of the $D_c * ln(r)$ curve.

In practice of course, we do not have an infinite number of data at our disposal. Therefore 
the left and right parts of the $ln(C(r))$ against $ln(r)$ diagram will not be reliable, 
so that the most informative slopes between those extremes in the diagram have to be identified.
When the size of the data space increases, it will reach the dimension where it is effectively 
possible to compute the correlation dimension (obviously, working in a too low dimensional 
space does not allow to estimate a large dimension). As this dimension is unknown, the experience 
must be carried out for increasing dimensions of the data space; when the required level is reached 
or exceeded, the estimated correlation dimensions will remain identical (i.e. the curves will be 
parallel, which can be observed on the curves as a saturation effect). 

Figure 1 shows the results obtained with the 4900 known values of the CATS time series. With 
respect to Grassberger-Proccacia's procedure described above, the data $x_t$ are now the 
$p$-dimensional regressors defined in (\ref{pred}). The figure shows a plot of $ln(C(r))$ against 
$ln(r)$ for increasing dimensions of the data space (i.e. increasing sizes $p$ of the regressors). 
The expected saturation effect \cite{GrassNPL} is clear for values of $ln(r)$ between 5 and 8; 
the correlation dimension seems to be around 1. 

\begin{figure}[!hbt]
\begin{center}
\scalebox{0.35}{\includegraphics*{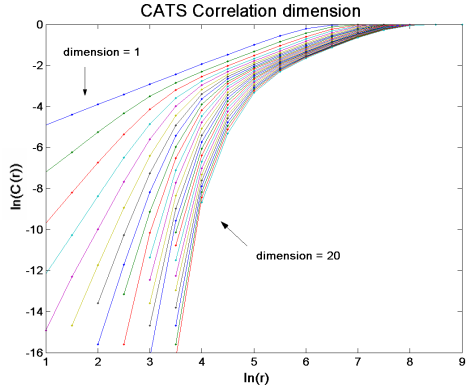}}
\end{center}
\caption{Estimation of the correlation dimension using the Grassberger-Procaccia procedure; 
log of the correlation integral $C(r)$ versus the log of the hypersphere radius $r$.}
\label{Grass_b_b}
\end{figure}

Another representation of the correlation dimension can be given in plotting the estimation 
of the correlation dimension as in Figure 2. A flat region can be seen around $ln(r) = 6$ to 7 
where the correlation dimension is again approximately one. In conclusion, according to Taken's 
theorem, any regressor for the CATS time series should be at most of size 3.

\begin{figure}[!hbt]
\begin{center}
\scalebox{0.35}{\includegraphics*{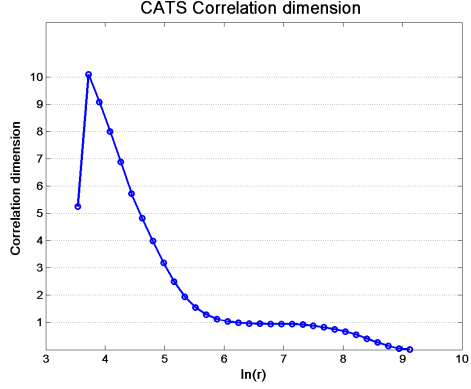}}
\end{center}
\caption{Correlation dimension obtained for various values of the hypersphere
radius $r$ (in log scale) for the CATS time series.}
\label{Grass_d_b}
\end{figure}

Note that the correlation dimension estimation is only a preliminary rough calculation, in order 
to get a first insight on the series. Indeed, because of the high correlation between successive 
values in the series (or in other words its 'smoothness'), it may happen that the correlation 
dimension estimation just catches this correlation, and not the dynamics of the series. 
In Figure \ref{Grass_b_b}, this could be seen in the form of two saturation effects in the slopes, 
as detailed above: one when the correlation between successive values is caught, the other one 
when the true dimensionality of the series is reached. According to the very low value (one) 
found for the correlation dimension in the CATS series, this risk does not have to be underestimated. 
Nevertheless, as no other reasonable value can be found, we will consider in the following that 
the value found for the correlation dimension is reliable, and a regressor of size 3 will thus 
be used.

The problem of a time series with a very low correlation dimension is that each value only 
depends on the few preceding ones. Any model built according to the above principles is 
therefore restricted to a very limited amount of information, and the prediction becomes hard 
and unstable. This is for example the case in financial time series prediction: the high sampling 
frequency of financial indexes makes them extremely smooth at short term. In such context, 
one usually models pre-processed series instead of the original ones; the preprocessing can consist 
in differences, returns, etc. Because of the similarities between the correlation dimension results 
on such financial series and on the CATS one, the same kind of pre-processing is developed here. 
In addition to the original series, two pre-processed ones will be used in the experiments: 
the series obtained by differences and by returns. The difference time series is obtained as:
\begin{equation}\label{serie_diff}
\displaystyle x_d(t) = x(t + 1) - x(t),
\end{equation}
while the return time series is computed as:
\begin{equation}\label{serie_returns}
\displaystyle x_r(t) = \frac{x(t + 1) - x(t)}{x(t)}.
\end{equation}

\begin{figure}[!hbt]
\begin{center}
\scalebox{0.35}{\includegraphics*{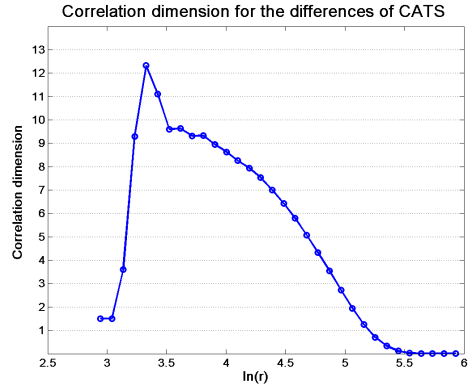}}
\end{center}
\caption{Correlation dimension obtained for various values of the hypersphere
radius $r$ (in log scale) for the difference time series.}
\label{Grass_d_d}
\end{figure}

\begin{figure}[!hbt]
\begin{center}
\scalebox{0.35}{\includegraphics*{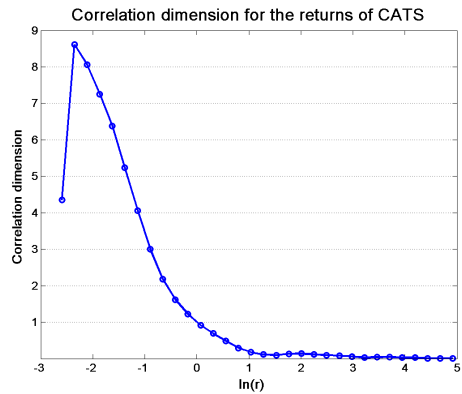}}
\end{center}
\caption{Correlation dimension obtained for various values of the hypersphere
radius $r$ (in log scale) for the return time series.}
\label{Grass_d_r}
\end{figure}

The correlation dimension of these two new time series can also be computed, using the same 
Grassberger-Procaccia procedure. The results are presented in Figures \ref{Grass_d_d} and 
\ref{Grass_d_r}. Unfortunately, in both cases, the results are not conclusive. Indeed 
there is no visible saturating slope in the $ln(C(r))$ against $ln(r)$ diagrams, contrarily 
to the saturation in Figure \ref{Grass_d_b}. The results found on the original series will thus 
be kept in the following as a rough estimation of the correlation dimension of the CATS time series.

\section{The double quantization forecasting method}\label{method}

\subsection{Self-organizing Maps}\label{SOM}

The Self-Organizing Map (SOM) is an unsupervised classification algorithm introduced in the 80's 
by Teuvo Kohonen \cite{Kohonen95}. Since their first description, Self-Organizing Maps have been 
applied in many different fields to solve various problems. Their theoretical properties are well 
established \cite{Cottrell98}, \cite{Cottrell97}. 

In a few words, a SOM map has a fixed number of units quantifying the data space. Those 
units, also called prototypes or centroids, are linked by predefined neighbourhood 
relationships that can be represented graphically through a 1- or 2-dimensional grid. 
After learning, the grid of prototypes has two properties. First, it defines a vector 
quantization of the input space, as any other vector quantization algorithm. Secondly, 
because the grid relationships are used in the learning algorithm itself, the grid 
representation has a topological property: two close inputs will be projected either 
on the same prototype or on two close ones in the grid. The Kohonen map can thus be seen 
as an unfolding procedure or as a nonlinear projection from the data space on a 1- or 
2-dimensional grid. The prototypes in a Kohonen map can also be seen as representatives 
of their associated class (the set of data nearer from a specific prototype than from 
any other one), turning the algorithm into a classification (or at least a clustering) 
tool. One of the main features of Kohonen maps is their ability to easily project data 
in a 2-dimensional representation, allowing intuitive interpretations.

\subsection{The Double Vector Quantization method (DVQ)}\label{DVQ}
Though the SOM is usually considered as a classification, feature extraction or recognition 
tool, there exist a few works where the SOM algorithm is used in time series prediction problems, 
as \cite{Walter90},\cite{Vesanto97}, \cite{Koskela98}, \cite{Cottrell96}, \cite{Cottrell}, 
\cite{Dablemont03}. In most of these situations however, the goal is to reach a reliable 
one-step-ahead prediction. In this work we are specifically looking for longer-term ones, 
and more precisely to 20 steps ahead prediction in the context of the CATS Competition.

A complete description of the DVQ method is given in \cite{Simon04}, together with a full proof 
of the method stability for long-term predictions. A brief description of the method is given 
here in the simple case of a scalar time series prediction. Full details for the vector case 
can be found in \cite{Simon04}. The goal of the method is to extract long-term information or 
trends of a time series. The method is based on the SOM algorithm used to characterize 
(or learn) the past of the series. Afterwards a forecasting step allows predicting future values.

\subsubsection{Characterization} 
According to the formulation of a nonlinear auto-regressive model (\ref{pred}), the method 
uses regressors of past values to predict the future evolution of a time series. Having 
at disposal a scalar time series of $n$ values, the correlation dimension $D_c$ is evaluated, 
leading to the choice of $p$-dimensional regressors. The $n$ known values of the time series 
are then transformed into $p$-dimensional regressors: 
\begin{equation}\label{regress}
\displaystyle x_t = \{x(t-p+1), ..., x(t-1), x(t)\},
\end{equation}
where $p \leq t \leq n$, and $x(t)$ is the original time series at our disposal. As one may 
expect $n-p+1$ such regressors are obtained from the original time series. 

The original regressors $x_t$ are then manipulated such that other regressors are created, 
according to:
\begin{equation}\label{deform}
\displaystyle y_t = x_{t+1}-x_t.
\end{equation}
The $y_t$ vectors are called the deformation regressors, or the deformations in short. 
By definition each deformation $y_t$ is associated to a single regressor $x_t$. Of course, 
$n-p$ deformations are obtained from a time series of $n$ values. 

At this stage of the method there exist two sets of regressors. The first one contains 
the $x_t$ regressors and is representative of the original space (of regressors). The space 
containing the $y_t$ deformations is representative of the deformation space. Those two sets 
of vectors, of the same dimension $p$, will be the data manipulated by the SOM maps.

Applying the SOM algorithm to each one of these two sets results in two sets of prototypes, 
denoted respectively $\bar{x}_i$, with $1 \leq i \leq n_1$, and $\bar{y}_j$, with 
$1 \leq j \leq n_2$. The classes associated to those prototypes are denoted respectively 
$c_i$ and $c'_j$. 

Characterizing the two time series through the quantization of the regressors and deformations 
is a static-only process. The dynamics of the past evolution of the series has to be modeled 
too. In fact, this is possible because the dynamics is implicitly recorded in the deformations. 
The issue is thus to build a representation of the existing relations between the original 
regressors and the deformations. For this purpose, a matrix $f(ij)$ is defined such that 
its $ij$ element, denoted $f_{ij}$, is obtained as:
\begin{equation}\label{freqs}
\displaystyle f_{ij} = \frac{\#\{x_t \in c_i \text{ and } y_t \in c'_j\}}{\#\{x_t \in c_i\}},
\end{equation}
with $1 \leq i \leq n_1$, $1 \leq j \leq n_2$. Intuitively the probability of having 
a certain deformation $j$ associated to a given regressor $i$ is approximated 
by the empirical frequencies (\ref{freqs}) measured on the data at disposal. Each row of the 
$f(ij)$ matrix ($1 \leq j \leq n_2$) in (\ref{freqs}) is in fact the conditional probability 
that $y_t$ belongs to $c'_j$ given the fact that $x_t$ belongs to $c_i$. Of course, 
elements $f_{ij}$ ($1 \leq j \leq n_2$) sum to one for each $i$.

\subsubsection{Forecasting}
Now that the past evolution of the time series has been modeled, predictions can be performed. 
Let us define the last known value $x(t)$ at time $t$, with corresponding regressor $x_t$. 
The prototype $\bar{x_k}$ closest to $x_t$ in the original space is searched. According to the 
conditional probability distribution defined by row $k$, a deformation prototype $\bar{y_l}$ 
is then chosen randomly among the $\bar{y_j}$, according to the $f_{kj}$ probability law. 
The prediction for instant $t+1$ is finally obtained according to relation (\ref{deform}):
\begin{equation}\label{pred+1}
\hat{x}_{t+1} = x_t + \bar{y}_l,
\end{equation}
where $\hat{x}_{t+1}$ is the estimate of $x_{t+1}$ given by the model. In fact $\hat{x}_{t+1}$ 
is a $p$-dimensional vector, and only one of its components corresponds to a prediction 
$\hat{x}(t+1)$ at time $t+1$; this value is thus extracted from the $\hat{x}_{t+1}$ vector 
and taken as the prediction. 

Once a one-step-ahead prediction (horizon $h$ = 1) is computed, the whole procedure can be 
repeated to obtain predictions for higher values of $h$. In practice, prediction $\hat{x}(t+1)$ 
is used to compute $\hat{x}_{t+2}$ through its corresponding regressor $\hat{x}_{t+1}$. 
$\hat{x}(t+2)$ is then extracted from $\hat{x}_{t+2}$, and so on up to horizon $h$. This 
recursive procedure is the standard way to obtain long-term forecasts from a one-step-ahead 
method. The whole procedure up to horizon $h$ is called a {\it simulation}.

\subsubsection{Comments} 
The goal of the DVQ method is to provide insights over the possible long-term evolution 
of a series, and not necessarily a single accurate prediction. The long-term (horizon $h$) 
simulations are then repeated using a Monte-Carlo procedure. The simulations distribution 
can be observed, and statistical information such as variance, confidence intervals, etc. 
can be determined too. The obtained long-term predictions have been proven to be stable 
\cite{Simon04}. 

Another important comment is that the method can easily be generalized to the prediction 
of vectors. With respect to the procedure described in the previous subsection, the only 
difference is that deformations (\ref{deform}) must be computed by differences of $d$-spaced 
values:
\begin{equation}\label{deform_d}
\displaystyle y_t = x_{t+d} - x_t,
\end{equation}
a direct generalization of the $d$ = 1 case in (\ref{deform}). Then, $d$ scalar values have 
to be extracted from the $\hat{x}_{t+d}$ vector, and so on. For example, two values could be 
extracted (corresponding to $\hat{x}(t+1)$ and $\hat{x}(t+2)$). In this case, repeating 
the procedure means to inject $\hat{x}(t+1)$ and $\hat{x}(t+2)$ to predict $\hat{x}(t+3)$ 
and $\hat{x}(t+4)$. More details about the vector case can be found in \cite{Simon04}.

A third comment concerns the numbers $n_1$ and $n_2$ of prototypes respectively in the regressor 
and deformation spaces. The major concern is that different values of $n_1$ ($n_2$) lead 
to different segmentations of the regressor and the deformation spaces which in turn lead 
to different models of the time series. The values of $n_1$ and $n_2$ should therefore 
be optimized, for example by resampling procedures (cross-validation, bootstrap) on a predefined 
error criterion, as the one-step-ahead error. 

Finally, since the only property of the SOM used here is the vector quantization, any other vector 
quantization method could have been chosen to implement the above procedure. The SOM maps have 
been chosen since they seem more efficient and faster compared to other VQ methods despite 
a limited complexity \cite{deBodt03}. Furthermore, they provide an intuitive and helpful graphical 
representation. Note that in practice any kind of SOM map could be used, but that one-dimensional 
maps, or strings, are preferred here.

\section{Methodological aspects of the double quantization for the CATS data set}\label{Methodo}
As mentioned in section \ref{DVQ} the goal of the DVQ method is to provide insights over the possible 
long-term evolution of a series, and not necessarily a single accurate prediction. In this section 
the methodology for the experiments will be described having in mind that the method has now to predict 
accurate values due to the competition context.

\subsection{Scalar and vector predictions}\label{scavect}
From section \ref{analyseCATS} we know that regressor $x_t$ for nonlinear models should 
contain at most 3 past values:
\begin{equation}\label{regress_3}
\displaystyle x_t = \{x(t-2), x(t-1), x(t)\}.
\end{equation}
As this expression has the same form as relation (\ref{regress}), it allows a direct application of 
the DVQ method to predict $x(t+1)$. This direct application of the method is an illustration of 
the scalar prediction with the DVQ method. 

The use of the correlation dimension to set the size of the regressor should however be taken 
with care when predicting vectors instead of scalars, as it will be the case for the CATS series. 
As an example, if one wants to predict a vector of $d = 2$ values, namely $\{\hat{x}(t+1), 
\hat{x}(t+2)\}$, the following two regressors should be respectively used:
\begin{equation}\label{pred_1_2}
\begin{array}{c}
\displaystyle \{x(t-2), x(t-1), x(t)\} \text{ to predict } \hat{x}(t+1),\\
\displaystyle \{x(t-1), x(t), x(t+1)\} \text{ to predict } \hat{x}(t+2).
\end{array}
\end{equation}
In order to use the DVQ method, it is suggested to merge the two regressors and use:
\begin{equation}\label{regress_pred_1_2+1}
\displaystyle \{x(t-2), x(t-1), x(t), x(t+1)\}
\end{equation}
to predict $\{\hat{x}(t+1), \hat{x}(t+2)\}$ together. Of course this is impossible since value 
$x(t+1)$ is unknown at time $t$. Therefore, the regressors that will be used in the following 
have the same size as in Eq. \ref{regress_pred_1_2+1} but they are shifted such that their last 
component now corresponds to $x(t)$:
\begin{equation}\label{regress_pred_1_2_fin}
\displaystyle \{x(t-3), x(t-2), x(t-1), x(t)\}.
\end{equation}
In this example, deformations $y_t$ are then computed according to:
\begin{equation}
\begin{tabular}{rcl}
$y_t$ & = & $x_{t+d}$ - $x_t$\\
 & = & $\{x(t-1), x(t), x(t+1), x(t+2)\}$ \\
 & & - $\{x(t-3), x(t-2), x(t-1), x(t)\}$.\\
\end{tabular}
\end{equation}
The same procedure can be extended for values of $d$ greater than 2.

To summarize, the DVQ method is directly applicable in the scalar case. Some care must be taken 
in the vector case: if vectors of $d$ values have to be predicted then the corresponding regressors 
have to be merged into a single vector which may only contain known values. Only then, the DVQ method 
in vector case can be applied.

\subsection{20 step ahead prediction strategies}\label{20steps}
As the CATS competition requires predicting block of 20 successive unknown values, several strategies 
can be used to reach this prediction horizon. The first one, called the {\it recursive strategy}, 
is probably the most common way to obtain a long term prediction. Predictions are obtained recursively 
until the final time horizon $h$; i. e. the one-step-ahead predictions are included one by one 
in the regressor to obtain the next one-step-ahead prediction. Formally, the last prediction $\hat{x}(t+k)$ 
is used to predict the next one $\hat{x}(t+k+1)$ as part of the regressor for $t+k+1$ ($1 \leq k \leq d-1$).

The second approach, called the {\it block strategy}, consists in predicting all the $h$ future values 
in one single vector. This is made possible by the use of a vector prediction method, as the DVQ one. 

A mixed approach is a {\it recursive-block strategy}, where blocks of intermediate size $d$ are 
predicted through a limited number of $h/d$ recursive steps (where $h$ is supposed to be a multiple 
of $d$ for simplicity).

\subsection{Number of prototypes}\label{n1n2}
As mentioned in section \ref{DVQ}, the numbers $n_1$ and $n_2$ of prototypes in respectively 
the regressor and the deformation spaces have to be fixed. A cross-validation procedure is therefore 
used. This cross-validation procedure mimics the competition problem. 

Fifteen new gaps of length 20 have been created randomly in the available data. As the true values 
for those 300 new missing values are known they can serve as validation set for models learned 
on the remaining values. Note that the random selection was slightly constrained: the new gaps 
could not overlap between them and with the existing ones. The whole validation procedure has been 
repeated 20 times, to average the dependencies to the choice of the 300 new gaps. 

To compare the different models that will be learned on the 20 various learning sets a mean square 
error $MSE$ validation criterion is used. This criterion is comparable to the one proposed 
in the CATS competition and is defined as:
\begin{equation}\label{MSE}
\displaystyle MSE = \frac{\displaystyle \sum_{y_t \in VS}^{}{(y_t - \hat{y}_t)^2}}{300},
\end{equation}
where $VS$ represents one of the 20 validation sets of 300 new missing values. The best model will 
be the one which has the lowest average $MSE$ over the 20 validation sets.

\subsection{Final predictions}\label{finalpreds}
Once the optimal $n_1$ and $n_2$ numbers are found, a new learning stage is done now using all 
available data (i.e. the 4900 data of the original CATS time series). 

To avoid problems due to the random initialization of the prototypes, several learning procedures 
are performed, and the best one is selected according to the validation sets, even if using 
the latter may lead to a small amount of overfitting (as the validation sets are now part 
of the new learning set). Simulations at final horizon are then repeated 100 times, and 
the mean is computed.

To refine this direct application of the DVQ method to the CATS prediction problem, some specific 
heuristics have been developed. The first one is to reverse the time series. Indeed, for the four 
blocks inside the series, the prediction can be performed from right to left (decreasing values of 
time) as well as it is performed from left to right (chronological order). For those four blocks 
of length 20, the CATS Competition is a missing value problem rather than a forecasting one. 
The DVQ method will therefore be applied in both directions.

The second heuristic considers a prediction horizon up to $h = 21$ (instead of $h = 20$). Indeed, 
as the 21st value is known for the four first gaps, predicting to horizon $h = 21$ allows a comparison 
between the true value and the 21st predicted one. As some error in long-term trends of the prediction 
is unavoidable, the comparison between these true and predicted values leads to an error that can be 
compensated at first order through a linear correction of the simulated predictions. This correction 
is performed such that the 21st predicted value is made equal to the true one.

The final predictions that were sent to the Competition were obtained by the DVQ method combined 
with the two heuristics explained above. More precisely, for the four first blocks of 20 missing 
values, a prediction to horizon $h = 21$ was performed in both directions, resulting from 100 DVQ 
simulations. As the true value for the 21st data is known in both directions for these four cases, 
a linear correction has been applied. The final predictions were then obtained as the average 
of the two linearly corrected predictions.

Of course, as the true 21st value is unknown for the fifth block of missing data, the above strategy 
can not be applied. The final predictions for this case were obtained as the mean of 100 DVQ simulations.

\begin{figure}[!hbt]
\begin{center}
\scalebox{0.6}{\includegraphics*{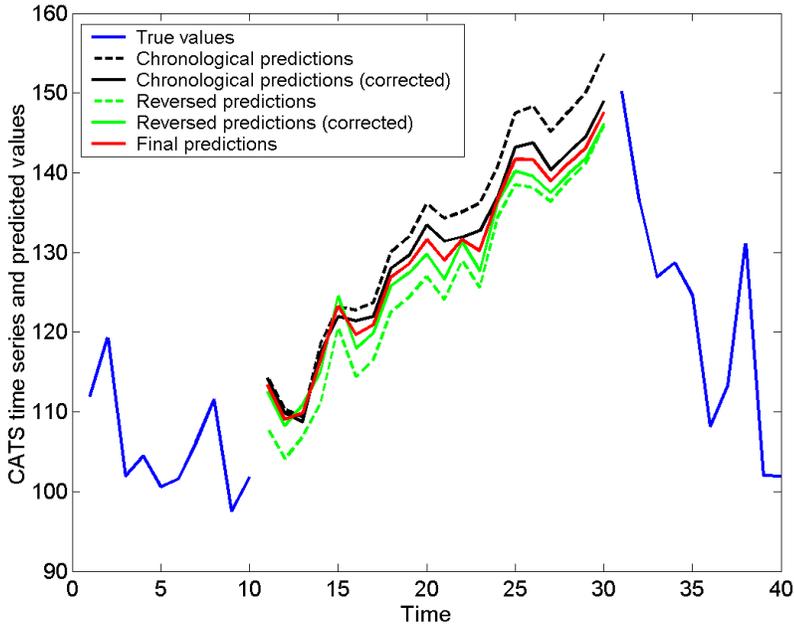}}
\end{center}
\caption{Corrections applied to obtain the final predictions using the two heuristics 
(first block of missing values; data 981-1000). See text for details.}
\label{heuristic}
\end{figure}

Figure \ref{heuristic} shows the various steps leading to the final predictions for the first block 
of missing values. The outmost curves are the mean of the 100 DVQ simulations (top of the figure 
is the chronological order, bottom is the reversed order). The two inner curves are the linearly 
corrected values using the comparison between the 21st true and predicted values. The fifth curve 
in the middle represents the final predictions, i.e. the mean of the two (linearly corrected) 
inner curves.

\section{Experimental results}\label{results}
According to the 'financial-like' behaviour of the CATS time series, as discussed in section 
\ref{analyseCATS}, three time series are considered in all our experiments: the initial CATS, 
the difference and the return time series. Furthermore, this 'financial-like' behaviour already 
suggests that a recursive strategy may behave poorly for a time horizon of 20 values. Consequently, 
in addition to the recursive strategy, where one-step-ahead predictions are repeated 20 times, 
a recursive-block strategy is used, with blocks of size 2, 5 and 10, and finally a block strategy 
is used with a bloc size equal to the time horizon $h = 20$. In the recursive-block strategy, 
the time horizon of 20 values corresponds to predict 10 blocks of size $d$ = 2, 4 blocks of size 5, 
etc.

For each one of the three time series, for each one of the block sizes, a cross-validation using 
the 20 validation sets has been performed as described in section \ref{n1n2}. For comparison purposes 
the new missing values in the 20 validation sets are identical in all experiments. Models with $n_1$ 
and $n_2$ both ranging from 5 to 100 by incremental steps of 5 are learned in each experiment. 
The $MSE$ criterion (\ref{MSE}) has been used to estimate the models generalization ability 
on the validation sets. 

Table \ref{resume_simus} gives a summary of the experiments. For each time series, for each block 
size, $n_1$ and $n_2$ corresponding to the best model in average are given, together with the average 
$MSE$. For the difference and return time series the $MSE$ is of course computed by first applying 
the inverse transformations on the predictions, in order to obtain MSE values that can be compared 
between the series. 

\begin{table}[!htp]
\begin{center}
\begin{tabular}{c|c|c|c|c}
Time series & \# step(s) ahead & $n_1$ & $n_2$ & $MSE$\\
\hline 
Initial & 1 & 90 & 5 & 1.66 $10^3$\\
 & 2 & 50 & 5 & 1.32 $10^3$\\
 & 5 & 25 & 5 & 1.36 $10^3$\\
 & 10 & 20 & 15 & 1.70 $10^4$\\
 & 20 & 65 & 10 & 2.98 $10^4$\\
\hline
Differences & 1 & 25 & 5 & 2.59 $10^3$\\
 & 2 & 90 & 5 & 1.90 $10^3$\\
 & 5 & 80 & 5 & 1.86 $10^4$\\
 & 10 & 55 & 5 & 4.67 $10^4$\\
 & 20 & 55 & 60 & 7.43 $10^4$\\
\hline
Returns & 1 & 10 & 5 & 3.39 $10^5$\\
 & 2 & 5 & 5 & 2.04 $10^5$\\
 & 5 & 55 & 5 & 1.83 $10^5$\\
 & 10 & 15 & 95 & 2.67 $10^{10}$\\
 & 20 & 45 & 55 & 4.01 $10^{10}$\\
\end{tabular}
\end{center}
\caption{Experiment summary: $n_1$ and $n_2$ for the best model in average over the 20 
cross-validations and corresponding $MSE$.}
\label{resume_simus}
\end{table}

From this table it seems clear that a model learned on the initial time series is adequate; none 
of the two preprocessing methods suggested in (\ref{serie_diff}) and (\ref{serie_returns}) reveals 
interesting. Indeed the MSE values obtained on the validation sets are in both cases larger 
than those obtained with the initial time series. 
	
\begin{figure}[!hbt]
\begin{center}
\scalebox{0.35}{\includegraphics*{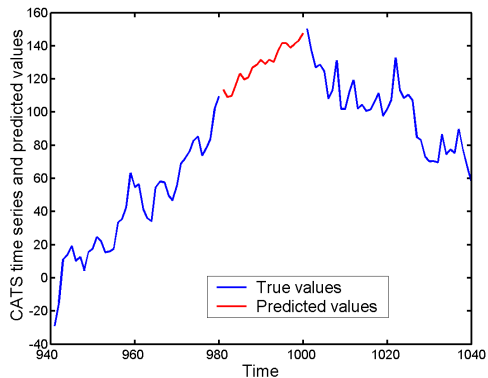}}
\end{center}
\caption{True values and final predictions, first gap.}
\label{hole1}
\end{figure}

\begin{figure}[!hbt]
\begin{center}
\scalebox{0.35}{\includegraphics*{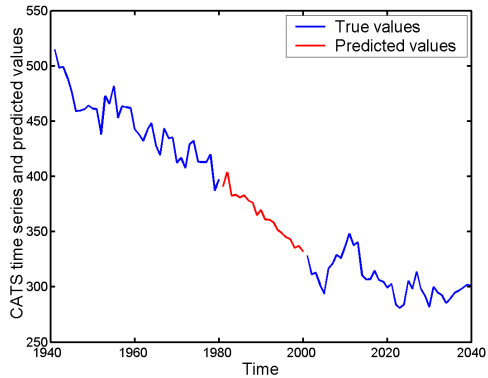}}
\end{center}
\caption{True values and final predictions, second gap.}
\label{hole2}
\end{figure}

\begin{figure}[!hbt]
\begin{center}
\scalebox{0.35}{\includegraphics*{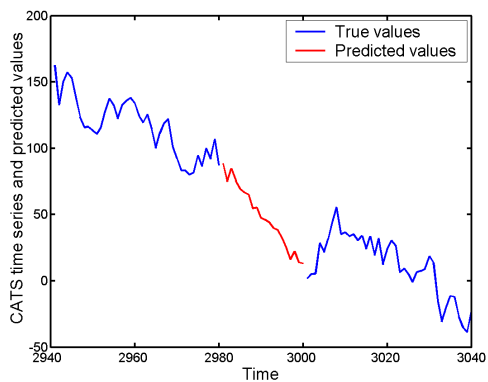}}
\end{center}
\caption{True values and final predictions, third gap.}
\label{hole3}
\end{figure}

\begin{figure}[!hbt]
\begin{center}
\scalebox{0.35}{\includegraphics*{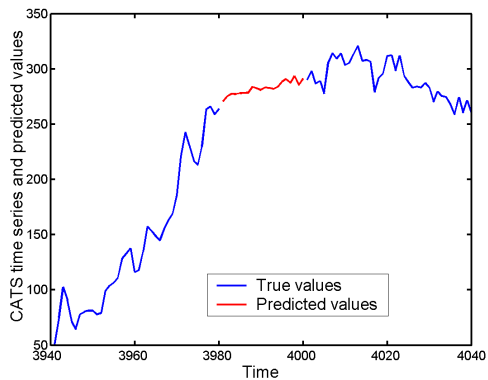}}
\end{center}
\caption{True values and final predictions, fourth gap.}
\label{hole4}
\end{figure}

\begin{figure}[!hbt]
\begin{center}
\scalebox{0.35}{\includegraphics*{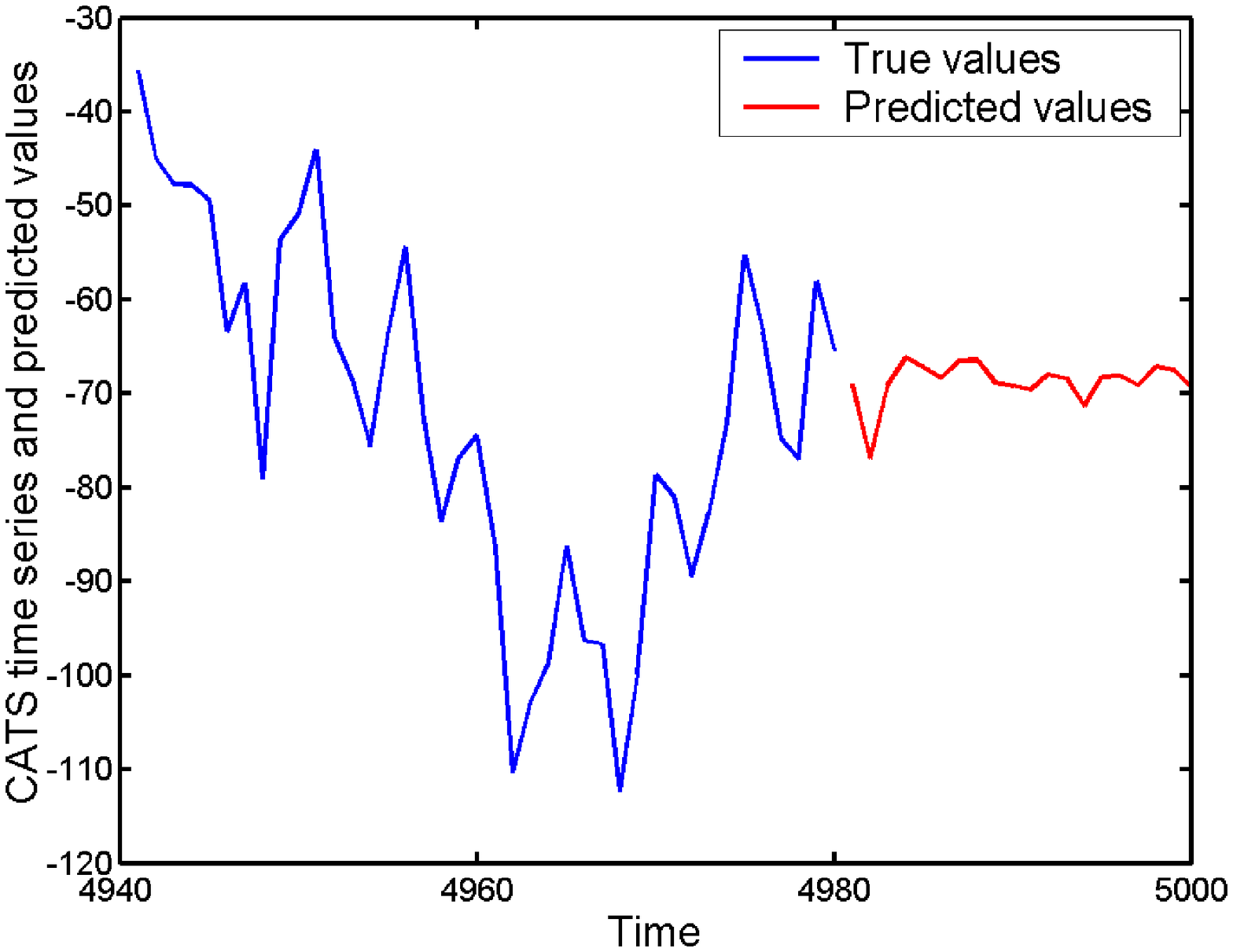}}
\end{center}
\caption{True values and final predictions, fifth gap.}
\label{hole5}
\end{figure}

Furthermore it is obvious that there exists a compromise between 20 repetitions of a one-step-ahead 
prediction (recursive strategy) and a single prediction of a vector containing the 20 next values 
(block strategy). This compromise seems to be somewhere between 10 predictions of blocks of 2 values 
and 4 predictions of blocks of 5 values. Nowadays the $MSE$ criterion is the lowest for blocks 
of size $d = 2$. The corresponding model, with 50 prototypes in the regressor space and 5 
in the deformation space, is selected to give the final prediction of the 100 missing values 
of the CATS Competition according to the heuristic described in section \ref{finalpreds}. 
Figures \ref{hole1} to \ref{hole5} show how the DVQ method predicts the missing values for the five 
20-values blocks of the CATS competition.

After the CATS Competition was closed, the results of the 17 classified competitors 
out of the 24 submissions were made available \cite{Lendasse04}. It can be seen in 
\cite{Lendasse04} that the DVQ method was ranked to the fourth position on the problem 
of predicting the first four gaps of 20 missing data. The value of the $E_2$ criterion in 
\cite{Lendasse04} (average sum of squares of errors on first 80 missing data) is 351. 
Besides the efficiency of the forecasting method, this is probably, and not surprisingly, 
a consequence of the fact that taking into account the first known value after each gap 
indeed improves the prediction accuracy. This information is not available for the fifth 
block of 20 missing data. The $E_1$ criterion in \cite{Lendasse04} (average sum of squares 
of errors on the 100 missing data) takes into account this very different problem. 
On this criterion, the DVQ method performs slightly worst with a result of 653, and 
is ranked in seventh position. 

\section{Conclusion}\label{ccl}
In this paper the results obtained with the double vector quantization method, based on the SOM maps, 
applied to the CATS data set are presented.

An analysis of the data shows some interesting aspects of the time series. Its correlation dimension 
seems to be as low as one. To take into account this particular aspect potentially limiting for 
nonlinear models other time series have been defined, i.e. the differences and the returns 
of the initial CATS series.

These three time series have been modeled using various sizes of prediction blocks corresponding 
to longer time horizons, in order to make the most of the vector prediction ability of the double 
vector quantization method. 

The number of units in the SOM maps has been discussed and selected using a cross-validation 
procedure on new gaps created randomly on the CATS data set. This procedure, together with 
the selected validation criterion, has been implemented to select the best model in average 
in conditions as close as possible to the Competition ones.

Heuristics specifically designed in the CATS Competition context are also described.

Finally, the predictions obtained using the best selected model and the heuristics are illustrated 
graphically.

\section*{Acknowledgements}
G. Simon is funded by the Belgian F.R.I.A., M. Verleysen is Senior Research Associate of the Belgian 
F.N.R.S. The scientific responsibility rests with the authors.

\begin{biography}
{\scalebox{0.7}{\includegraphics[clip,keepaspectratio]{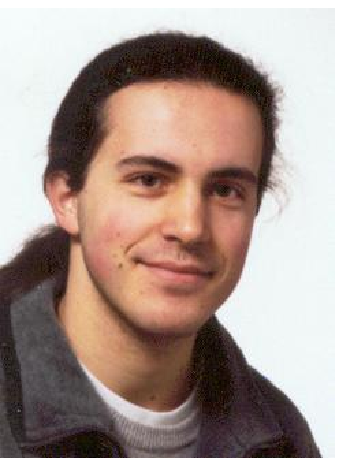}}}
{Geoffroy Simon}
{was born in 1978 in Dinant, Belgium. He received the M.Sc. degree in Computer Sciences 
in 2002 from the Facult\'es Universitaires Notre Dame de la Paix (Namur, Belgium). 
He is now working as Ph. D. student at the Microelectronic Laboratory of the Electrical 
Engineering Departement of the Universit\'e catholique de Louvain (UCL). His research topics 
cover nonlinear time series analysis, artificial neural networks, nonlinear statistics and 
self-organization applied to time-series forecasting problems. He is currently member 
of the UCL Machine Learning Group. His work is funded by a grant from the Belgian F.R.I.A. 
(Fonds pour la formation à la Recherche dans l'Industrie et dans l'Agriculture).}
\end{biography}

\begin{biography}
{\scalebox{0.38}{\includegraphics[clip,keepaspectratio]{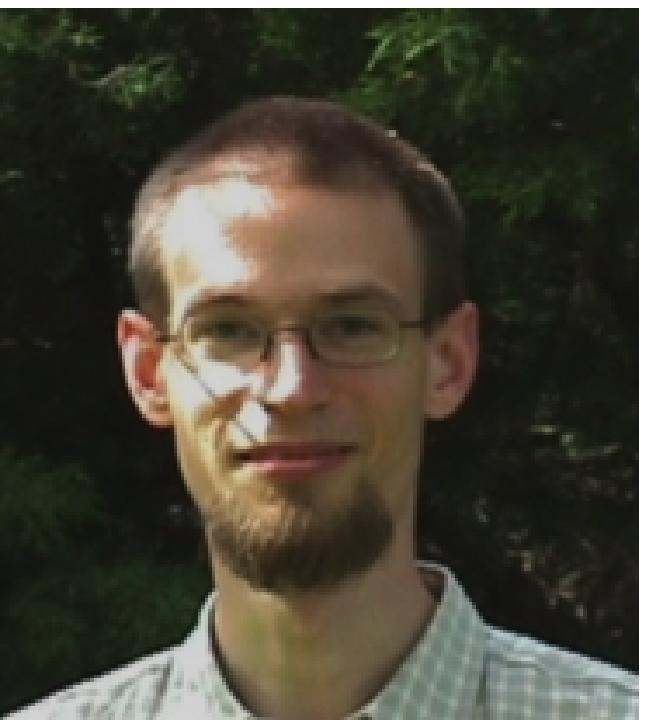}}}
{John Aldo Lee}
was born in 1976 in Brussels, Belgium. He received the M.Sc. degree in Applied Sciences 
(Computer Engineering) in 1999 and the Ph.D. degree in Applied Sciences (Machine Learning) 
in 2003, both from the Universit\'e catholique de Louvain (UCL, Belgium). His main interests 
are nonlinear dimensionality reduction, intrinsic dimensionality estimation, independent 
component analysis, clustering and vector quantization. He is a former member of the UCL 
Machine Learning Group and is now a Scientific Research Worker of the Belgian F.N.R.S. 
(Fonds National de la Recherche Scientifique). His current work aims at developing specific 
image enhancement techniques for Positron Emission Tomography in the Molecular Imaging and 
Experimental Radiotherapy department of the Saint-Luc University Hospital (Belgium).
\end{biography}

\begin{biography}
{\scalebox{0.61}{\includegraphics[clip,keepaspectratio]{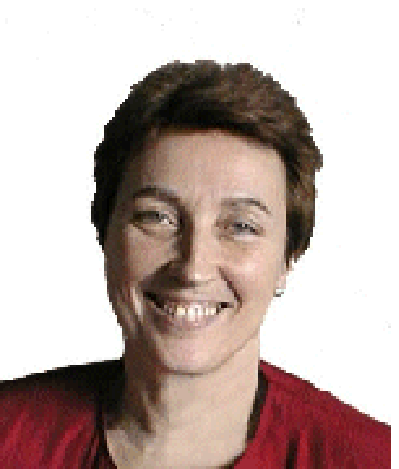}}}
{Marie Cottrell}
was born in 1943 in B\'ethune, France. She was a student at the Ecole Normale Sup\'erieure 
de S\`evres, and received the Agr\'egation de Math\'ematiques degree in 1964 (with 8th place), 
and the Th\`ese d'Etat (Mod\'elisations de r\'eseaux de neurones par des cha\^ines de Markov et 
autres applications) in 1988. From 1964 to 1967, she was a High School Teacher. From 1967 to 1988, 
she was successively an Assistant and an Assistant Professor at the University of Paris and 
at the University of Paris-Sud (Orsay), except from 1970 to 1973, on which she was a Professor 
at the University of Havana, Cuba. From 1989, she is a full Professor at the University Paris 1 
- Panth\'eon-Sorbonne. Her research interests include stochastic algorithms, large deviation theory, 
biomathematics, data analysis and statistics. Since 1986, her main works deal with artificial and 
biological neural networks, Kohonen maps and their applications in data analysis. She is the author 
of about 80 publications in this field. She is in charge of the SAMOS research group at the University 
Paris 1. She is regularly solicited as referee or international conference program committee member. 
In 2005, she organized the Fifth WSOM Conference in Paris.
\end{biography}

\begin{biography}
{\scalebox{0.45}{\includegraphics[clip,keepaspectratio]{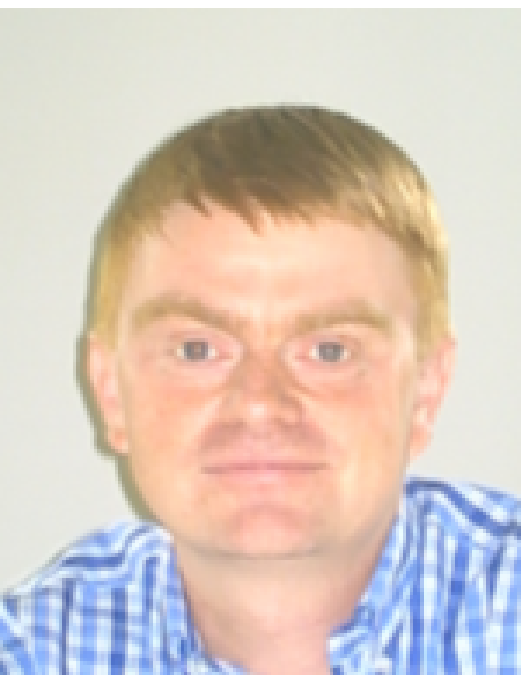}}}
{Michel Verleysen}
was born in 1965 in Belgium. He received the M.S. and Ph.D. degrees in electrical engineering 
from the Universit\'e catholique de Louvain (Belgium) in 1987 and 1992, respectively. He was an 
Invited Professor at the Swiss E.P.F.L. (Ecole Polytechnique F\'ed\'erale de Lausanne, Switzerland) 
in 1992, at the Universit\'e d'Evry Val d'Essonne (France) in 2001, and at the Universit\'e Paris 1 
- Panth\'eon-Sorbonne in 2002, 2003 and 2004.  He is now a Research Director of the Belgian F.N.R.S. 
(Fonds National de la Recherche Scientique) and Lecturer at the Universit\'e catholique de Louvain. 
He is editor-in-chief of the Neural Processing Letters journal, chairman of the annual ESANN 
conference (European Symposium on Artificial Neural Networks), associate editor of the IEEE Trans. 
on Neural Networks journal, and member of the editorial board and program committee of several 
journals and conferences on neural networks and learning. He is author or co-author of about 200 
scientific papers in international journals and books or communications to conferences with 
reviewing committee. He is the co-author of the scientific popularization book on artificial 
neural networks in the series "Que Sais-Je?", in French. His research interests include machine 
learning, artificial neural networks, self-organization, time-series forecasting, nonlinear statistics, 
adaptive signal processing, and high-dimensional data analysis.\end{biography}

\end{document}